\documentclass[12pt]{article} 
\usepackage{amsmath,amsxtra,amssymb,latexsym, amscd,amsthm}

\advance\voffset-1truecm\relax
\advance\hoffset-1.0truecm\relax
\font\titbf=cmbx10 scaled \magstep2

\font\tac=cmcsc10 scaled \magstep1

\numberwithin{equation}{section}

\newtheorem{theorem}{Theorem}

\def\C{\Bbb C}
\def\Z{\Bbb Z}

\begin{document}
$\qquad$

\centerline{\titbf  A UNIQUENESS THEOREM FOR MEROMORPHIC }
\medskip

\centerline{\titbf  MAPPINGS WITH TWO FAMILIES OF HYPERPLANES }
\medskip
\centerline{\tac {\bf Gerd Dethloff and Si Duc Quang and Tran Van Tan}}

\vskip0.15cm

\begin{abstract}
\noindent  In this paper, we extend the uniqueness theorem for meromorphic mappings to the case where the family of hyperplanes depends on the meromorphic mapping and where the meromorphic mappings may be degenerate.
\end{abstract}

\section{Introduction}
The uniqueness problem of meromorphic mappings under a condition
on the inverse images of divisors was first studied by Nevanlinna \cite{Ne}. He
showed that for two nonconstant meromorphic functions $f$ and $g$ on the
complex plane $\mathbb{C}$, if they have the same inverse images for five
distinct values, then $f\equiv g.$  In 1975,  Fujimoto \cite{F1} generalized
 Nevanlinna's result to the case of meromorphic mappings of $\mathbb{C}^{m}$
into $\mathbb{C}P^{n}$. He showed  that for two linearly nondegenerate meromorphic
mappings $f$ and $g$ of $\mathbb{C}^{m}$ into $\mathbb{C}P^{n}$, if they
have the same inverse images counted with multiplicities for $(3n+2)$
hyperplanes in general position in $\mathbb{C}P^{n},$ then $f\equiv g$. 

In 1983, Smiley \cite{Sm} showed that
\begin{theorem}\label{thSmiley}
 Let $f, g$ be linearly nondegenerate meromorphic mappings  of $\mathbb{C}^m$ into $\mathbb{C} P^n.$  Let $\{H_j\}_{j=1}^q$  ($q\geq 3n+2)$ be hyperplanes in $\mathbb{C} P^n$ in general position. Assume that

$\text{a)}\quad$   $f^{-1}(H_j)= g^{-1}(H_j)\ ,\quad \text{for all}\ \ 1\leq j\leq  q$ (as sets),

$\text{b)}\quad$   $\dim \big(f^{-1}(H_i)\cap f^{-1}(H_j)\big)\le m-2\ \ \text{for all}\  1\le i<j\le q\,$,

$\text{c)}\quad$  $f=g$ on \ \
$\bigcup_{j=1}^{q}f^{-1}(H_j)\,$.

\noindent Then $f\equiv g.$
\end{theorem}

\hbox to 5cm {\hrulefill }

 \noindent {\small Mathematics Subject
Classification 2000: Primary 32H30; Secondary 32H04, 30D35.}

\noindent {\small Key words :  Meromorphic mappings, uniqueness theorems.}

\noindent {\small  Acknowledgement: This work was done during a stay of the third named author at the Institut des Hautes \'Etudes Scientifiques, France. He
wishes to express his gratitude to this institute.}

\newpage

\noindent In 2006 Thai-Quang \cite{TQ} generalized this result of Smiley to the case where $q\geq 3n+1$ and $n\geq 2.$ In 2009, Dethloff-Tan \cite{DT} showed that for every nonnegative integer $c$ there exists a positive integer $N(c)$ depending only on   $c$ such that  Theorem \ref{thSmiley} remains valid if $q\geq (3n+2 -c)$ and $n\geq N(c)$. They  also showed that the coefficient of $n$ in the formula of $q$ can be replaced by a number which is smaller than 3 for all $n>>0.$ Furthermore, they established a uniqueness theorem for the case of $2n+3$ hyperplanes and multiplicities are truncated by $n.$ At the same time, they strongly generalized many uniqueness theorems of previous authors such as Fujimoto \cite{F2}, Ji \cite{Ji} and  Stoll \cite{St}.  Recently, by using again the technique of Dethloff-Tan \cite{DT}, Chen-Yan \cite{CY} showed that the assumption ``multiplicities are truncated by $n$" in the result of Dethloff-Tan can be replaced by  ``multiplicities are 
truncated by $1$".  In \cite{Q}, Quang examined the uniqueness problem for the case of $2n+2$ hyperplanes.   

We would like to note that so far, all results on the uniqueness problem have still been  restricted to the case where meromorphic mappings are sharing a common  family of hyperplanes. The purpose of this paper is to introduce a uniqueness theorem for the case where the family of hyperplanes depends on the meromorphic mapping. We also will allow that  the meromorphic mappings may be degenerate. For this purpose we introduce some new techniques which can also be used to obtain  simpler proofs for many other uniqueness theorems.\\

 We shall prove the following uniqueness theorem:

\begin{theorem} \label{TQT} 
Let $f, g$ be nonconstant  meromorphic mappings  of $\mathbb{C}^m$ into $\mathbb{C} P^n.$  Let $\{H_j\}_{j=1}^q$ and $\{L_j\}_{j=1}^q$ $(q>2n+2)$  be families of hyperplanes in $\mathbb{C} P^n$ in general position. Assume that 

$\text{a)}\quad$   $f^{-1}(H_j)= g^{-1}(L_j)\ \quad \text{for all}\ \ 1\leq j\leq  q\,,$

$\text{b)}\quad$   $\dim \big(f^{-1}(H_i)\cap f^{-1}(H_j)\big)\le m-2\ \ \text{for all}\  1\le i<j\le q$\,, 

$\text{c)}\quad$  $\frac{(f,H_i)}{(g,L_i)}=\frac{(f,H_j)}{(g,L_j)}$ on 
$\bigcup_{k=1\atop }^qf^{-1}(H_k)\setminus\big(f^{-1}(H_i)\cup f^{-1}(H_j)\big)$ for all $1\leq i<j\leq q\,.$

Then the following assertions hold :

$i)\quad\dim\langle \text{Im}f\rangle=\dim\langle \text{Im}g\rangle\stackrel{\text{Def.}}{=:}p,$

\noindent where for a subset $X\subset \C P^n,$ we denote by $\langle X\rangle$ the smallest projective subspace of $\C P^n$ containing $X.$

$ii) \quad$  If
$$(*)\quad q>\frac{2n+3-p+\sqrt{(2n+3-p)^2+8(p-1)(2n-p+1)}}{2}\;(\geq 2n+2),\,$$
then $$\frac{(f,H_1)}{(g,L_1)}\equiv\cdots\equiv\frac{(f,H_q)}{(g,L_q)}\,.$$ Furthermore,  there exists a linear projective transformation $\mathcal L$ of $\C P^n$ into itself  such that $\mathcal L(f)\equiv g$ and $\mathcal L(H_j\cap \langle \text{Im}f\rangle)=L_j\cap \mathcal L(\langle \text{Im}f\rangle)$ for all $j\in\{1,\dots,q\}.$ 

\end{theorem}

\noindent {\bf Remark.}
1.) In Theorem \ref{TQT} condition $c)$ is well defined since, by condition $a)$, 
$\frac{(f,H_i)}{(g,L_i)}$ is a (nonvanishing) holomorphic function outside $f^{-1}(H_i)$.

2.) The condition $(*)$ is satisfied in the following cases:

     $+) \quad q\geq 2n+3$ and $p\in\{1,2,n-1,n\}, n\in\Z^+.$

 $+) \quad q\geq2n+p+1$ and $ p\in\{2, 3, \dots,n\}, n\in \Z^+.$

3.) If there exists a subset $\{j_0,\dots,j_n\}\subset\{1,\dots,q\}$ such that $H_{j_i}\equiv L_{j_i}$ for all $i\in\{0,\dots,n\},$ then the proof of  Theorem \ref{TQT} implies that $f\equiv g.$ 

4.) For the special case where $f, g$ are linearly nondegenerate 
(i.e. $p=n$) and $H_j\equiv L_j\,$, from Theorem \ref{TQT} we get again the results of Dethloff-Tan \cite{DT} and Chen-Yan \cite{CY}.

\section{Preliminaries}

We set $\Vert z \Vert := \big(|z_1|^2 + \cdots + |z_m|^2\big)^{1/2}$
for $z = (z_1, \dots, z_m) \in \C^m$ and define
$$B(r) := \big\{z \in \C^m : \Vert z \Vert < r\big\},\qquad
S(r) := \big\{z \in \C^m : \Vert z \Vert = r\big\}
$$ 
for all $0 < r < \infty$.
Define
\begin{align*}
d^c &:= \frac{\sqrt{-1}}{4\pi}(\overline \partial - \partial),\quad
\upsilon := \big(dd^c\Vert z \Vert^2 \big)^{m-1}\quad \\
\sigma &:= d^c \text{log}\Vert z\Vert^2 \land 
\big(dd^c \text{log}\Vert z\Vert^2\big)^{m-1}.
\end{align*}

Let $F$ be a nonzero holomorphic function on $\C^m$.
For each $a \in \C^m$, expanding $F$ as $F = \sum P_i(z-a)$
with homogeneous polynomials $P_i$ of degree $i$ around $a$,
we define
$$ \nu_F(a) := \min \big\{ i : P_i \not\equiv 0 \big\}.$$
Let  $\varphi$ be a nonzero meromorphic function on $\C^m$.
We define the zero divisor $\nu_\varphi$ as follows: For each $z \in 
\C^m$, we choose nonzero holomorphic functions $F$ and $G$ on a 
neighborhood $U$ of $z$ such that $\varphi = {F}/{G}$ on
$U$ and $\text{dim}\,\big(F^{-1}(0) \cap G^{-1}(0)\big) \leqslant m-2$.
Then we put $\nu_\varphi(z) := \nu_F(z)$.

\noindent Let $\nu$ be a divisor in $\mathbb{C}^m$ and $k$ be positive integer or $+\infty $. Set \ $
|\nu|:=\overline{\big\{z:\ \nu(z)\neq 0
\big\}}$  and $\nu^{[k]}(z) := \min \{ \nu (z), k\}.$ 

The truncated counting function of $\nu$ is defined by
$$ N^{[k]}(r, \nu) := \int\limits_1^r 
\frac{n^{[k]}(t)}{t^{2m-1}} dt \quad (1 < r < + \infty), $$
where
\begin{align*}
n^{[k]}(t) = \begin{cases}
\displaystyle{\int\limits_{|\nu| \cap B(t)}}
\nu^{[k]}\cdot \upsilon \ &\text{for}\ m \geqslant 2,\\
\sum\limits_{|z| \leqslant t} \nu^{[k]}(z)
&\text{for}\ m = 1.\end{cases}
\end{align*}
We simply write \ $N(r, \nu)$ for $N^{[+\infty]}(r,\nu)$.

\noindent For a nonzero meromorphic function $\varphi$ on $\mathbb{C}^{m},$ we set \quad $N_{\varphi}^{[k]}(r):=N^{[k]}(r, \nu_\varphi)$ and $ N_{\varphi}(r):=N^{[+\infty]}(r, \nu_\varphi).$ 
We have the following Jensen's formula:
$$ N_\varphi (r) - N_{\frac{1}{\varphi}}(r) =
\int\limits_{S(r)} \text{log}|\varphi| \sigma 
- \int\limits_{S(1)} \text{log}|\varphi| \sigma .$$

Let $f : \C^m \longrightarrow \C P^n$ be a meromorphic mapping.
For an arbitrary fixed homogeneous coordinate system 
 $(w_0 : \cdots : w_n)$ in
$\C P^n$, we take a reduced representation $f = (f_0 : \cdots : f_n)$,
which means that each $f_i$ is a holomorphic function on $\C^m$
and $f(z) = (f_0(z) : \cdots : f_n(z))$ outside the analytic set
$\{ f_0 = \cdots = f_n = 0\}$ of codimension $\geqslant 2$.
Set $\Vert f \Vert = \big(|f_0|^2 + \cdots + |f_n|^2\big)^{1/2}$.
The characteristic function $T_f(r)$ of $f$ is defined by
$$ T_f(r) := \int\limits_{S(r)} \text{log} \Vert f \Vert \sigma -
\int\limits_{S(1)} \text{log} \Vert f \Vert \sigma , \quad
1 < r < + \infty . $$   
For a meromorphic function $\varphi$ on $\C^m$, the characteristic
function $T_\varphi(r)$ of $\varphi$ is defined by considering
$\varphi$ as a meromorphic mapping of $\C^m$ into $\C P^1$.

We state the First Main Theorem and the Second Main Theorem in Value Distribution 
Theory:
For a hyperplane $H : a_0 w_0 + \cdots + a_n w_n = 0$ in $\C P^n$
with Im$f \not\subseteq H$, we put $(f,H) = a_0 f_0 + \cdots 
+ a_n f_n$, where $(f_0 : \cdots : f_n)$  is a reduced
representation of $f$.\\
\medskip
\noindent
{\bf First Main Theorem.} 
 {\it Let $f$ be a meromorphic mapping of
$\C^m$ into $\C P^n$, and $H$ be a hyperplane in $\C P^n$ such that $(f,H) \not\equiv 0$. Then}
$$ N_{(f,H)}(r) \leqslant T_f(r) + O(1) \quad \text{\it for all}\ r > 
1.$$ Let  $n,N,q$ be positive integers with $q\geq 2N-n+1$ and $N\geq n.$  We say that hyperplanes $H_1,\dots,H_q$  in $\C P^n$  are in $N$-subgeneral position if $\cap_{i=0}^N H_{j_i}=\varnothing$ for every subset $\{j_0,\dots,j_N\}\subset\{1,\dots,q\}.$

\medskip
\noindent
{\bf Cartan-Nochka Second Main Theorem (\cite{No}, Theorem 3.1).} {\it Let $f$ be a linearly nondegenerate
meromorphic mapping of $\C^m$ into $\C P^n$ and
$H_1, \dots, H_q$  hyperplanes in $\C P^n$
in $N$-subgeneral position $(q\geq 2N-n+1).$ Then
$$ (q-2N+n-1) T_f(r) \leqslant \sum_{j=1}^q N_{(f,H_j)}^{[n]}(r) +
o\big(T_f(r)\big) $$ 
for all $r$ except for a subset $E$ of $(1, +\infty)$ of finite 
Lebesgue measure.}

\section{Proof of Theorem 3 }
We first remark that $f^{-1}(H_j)=g^{-1}(L_j)\ne\C P^n$ for all $j\in\{1,\dots,q\},$  and that therefore $\{H_j\cap  \langle\text{Im}f\rangle\}_{j=1}^q$ (respectively $\{L_j\cap \langle\text{Im}g\rangle\}_{j=1}^q)$ are hyperplanes in $\langle\text{Im}f\rangle$ (respectively $\langle\text{Im}g\rangle$) in $n-$subgeneral position: Indeed, otherwise there exists $t\in\{1,\dots,q\}$ such that $f^{-1}(H_t)=\C P^n.$ Then by the assumption $b)$ we have dim$f^{-1}(H_j)\leq m-2$ for all $j\in\{1,\dots,q\}\setminus\{t\}.$ Therefore, $f^{-1}(H_j)=\varnothing$ for all $j\in\{1,\dots,q\}\setminus\{t\}.$ Then $ \langle\text{Im}f\rangle\not\subset H_j$ for all $j\in\{1,\dots,q\}\setminus\{t\}.$ Thus, $ \{H_j\cap \langle\text{Im}f\rangle\}_{j=1\atop j\ne t}^q$ are hyperplanes in $\langle\text{Im}f\rangle$ in $n$-subgeneral position.

\noindent  By  the Cartan-Nochka Second Main Theorem, we have
\begin{align*}
(q-2n+\dim \langle\text{Im}f\rangle-2)T_f(r)\leq \sum_{j=1\atop j\ne t}^q N_{(f, H_j)}^{[\dim \langle\text{Im}f\rangle]}(r)+o(T_f(r)) =o(T_f(r)).
\end{align*}
This is a contradiction to the fact that $q > 2n+2$.

Since $\{H_j\}_{j=1}^{n+1}$ and $\{L_j\}_{j=1}^{n+1}$ are families of hyperplanes in general position, $\tilde{f}:=\big((f,H_1):\cdots :(f,H_{n+1})\big)$ and $\tilde{g}:=\big((g,L_1):\cdots:(g,L_{n+1})\big)$ are reduced representations of meromorphic mappings $\tilde{f}$ and $\tilde{g}$ respectively of $\C^m$ into $\C P^n.$ Furthermore, $\dim\langle\text{Im}f\rangle=\dim\langle\text{Im}\tilde{f}\rangle,$  $\dim\langle\text{Im}g\rangle=\dim\langle\text{Im}\tilde{g}\rangle,$  $T_{\tilde{f}}(r)=T_f(r)+O(1)$ and  $T_{\tilde{g}}(r)=T_g(r)+O(1).$

 By assumptions $a)$ and $c)$ we that 
\begin{align}\label{p1}
\tilde{f}=\tilde{g}\quad \text{on}\;\cup_{j=1}^qf^{-1}(H_j).
\end{align}
We now prove that
\begin{align}\label{p2}
\dim \langle\text{Im}f\rangle=\dim\langle\text{Im}g\rangle\stackrel{\text{Def.}}{=}p.
\end{align}
This is  equivalent to prove that $\dim \langle\text{Im}\tilde{f}\rangle=\dim\langle\text{Im}\tilde{g}\rangle .$ 
Therefore, it suffices to show  that for any hyperplane $H$ in $\C P^n$ then
\begin{align*}
(H,\tilde{f})\equiv 0\quad\text{if and only if}\quad (H,\tilde{g})\equiv 0.
\end{align*}
Suppose that the above assertion does not hold. Without loss of the generality, we may assume that there exists a hyperplane $H$ such that $(H, \tilde{f})\not\equiv 0$ and $(H,\tilde{g})\equiv 0.$  Then by (\ref{p1}) we have 
\begin{align}\label{p3}
(\tilde{f},H)=0\quad\text{on}\;\cup_{j=1}^qf^{-1}(H_j).
\end{align}
By (\ref{p3}) and by the First Main Theorem and the Cartan-Nochka Second Main Theorem we have
\begin{align*}
(q-2n+\dim \langle\text{Im}f\rangle-1)T_{f}(r)+O(1)
&\leq\sum_{j=1}^qN_{(f,H_j)}^{[\dim \langle\text{Im}f\rangle]}(r)+o(T_f(r))\\
&\leq \dim \langle\text{Im}f\rangle\sum_{j=1}^qN_{(f,H_j)}^{[1]}(r)+o(T_{f}(r))\\
&\stackrel{(\ref{p3})}{\leq}\dim \langle\text{Im}f\rangle N_{(\tilde{f},H)}(r)+o(T_{f}(r))\\
&\leq \dim \langle\text{Im}f\rangle T_{\tilde{f}}(r)+o(T_f(r))\\
&=\dim \langle\text{Im}f\rangle T_{f}(r)+o(T_f(r)).
\end{align*}
This is a contradiction to the fact that $q > 2n+2$.
We complete the proof of  (\ref{p2}). 

Now we prove that 
\begin{align}\label{a0}
\frac{(f,H_1)}{(g,L_1)}\equiv\cdots\equiv\frac{(f,H_q)}{(g,L_q)}.
\end{align}
We distinguish the following two cases:

{\bf Case 1:} There exists a subset $J:=\{j_0,\dots,j_{n}\}\subset\{1,\dots,q\}$ such that $$\frac{(f,H_{j_0})}{(g,L_{j_{0}})}\equiv\cdots\equiv\frac{(f,H_{j_{n}})}{(g,L_{j_{n}})}\stackrel{\text{Def.}}{\equiv}u\,.$$ 

We have Pole$(u) \cup $Zero$(u)\subset f^{-1}(H_{j_0})\cap f^{-1}(H_{j_1})$, which is an analytic set of codimension at least 2 by assumption $b)$.  Hence,  $\text{Pole}(u) \cup  \text{Zero} (u)=\varnothing.$

Since $H_{j_0},...,H_{j_n}$ are hyperplanes in general position, 
$F:=\big((f,H_{j_0}):\cdots:(f,H_{j_{n}})\big)$
is the reduced representation of a meromorphic mapping $F$ of $\C^m$ into $\C P^n$. Still by the same reason $T_F(r)=T_f(r)+O(1).$

Suppose that (\ref{a0}) does not hold. Then, there exists  $i_0\in\{1,\dots,q\}\setminus\{j_0,\dots,j_n\}$ such that 
\begin{align}\label{h1}
\frac{(f,H_{i_0})}{(g,L_{i_0})}\not\equiv u.
\end{align} 
Since the families $\{H_j\}_{j=1}^q$ and $\{L_j\}_{j=1}^q$ are in general position,  there exist hyperplanes $H^{i_0}: a_0\omega_0+\cdots+a_n\omega_n=0, \;L^{i_0}:b_0\omega_0+\cdots+b_n\omega_n=0$  in $\C P^n$ such that
$(f,H_{i_0})\equiv (F,H^{i_0}),$ and $(g,L_{i_0})\equiv b_0(g,L_{j_0})+\cdots +b_n(g,L_{j_n})\equiv\frac{(F,L^{i_0})}{u}.$
Therefore, by (\ref{h1})  we have
\begin{align*}
\frac{(F,H^{i_0})}{(F,L^{i_0})}\equiv \frac{(f,H_{i_0})}{u(g,L_{i_0})}\not\equiv 1.
\end{align*}
By assumption $c)$ and since $\text{Pole}(u) \cup  \text{Zero} (u)=\varnothing$, we have 
$u=\frac{(f,H_{j_0})}{(g,L_{j_0})}=\frac{(f,H_{i_0})}{(g,L_{i_0})}=u\frac{(F,H^{i_0})}{(F,L^{i_0})}$ on  $\big(\bigcup_{k=1}^qf^{-1}(H_k)\big)\setminus\big(f^{-1}(H_{i_0})\cup f^{-1}(H_{j_0})\big)$ and  $u=\frac{(f,H_{j_1})}{(g,L_{j_1})}=\frac{(f,H_{i_0})}{(g,L_{i_0})}=u\frac{(F,H^{i_0})}{(F,L^{i_0})}$ on  $\big(\bigcup_{k=1}^qf^{-1}(H_k)\big)\setminus\big(f^{-1}(H_{i_0})\cup f^{-1}(H_{j_1})\big).$ Then $\frac{(F,H^{i_0})}{(F,L^{i_0})}=1$ on $\big(\bigcup_{k=1}^qf^{-1}(H_k)\big)\setminus f^{-1}(H_{i_0})$.

\noindent Therefore,  
\begin{align*}
\sum_{k=1,k\ne i_0}^qN_{(f,H_k)}^{[1]}(r)&\leq N_{\frac{(F,H^{i_0})}{(F,L^{i_0})}-1}(r)\\
&\leq T_{\frac{(F,H^{i_0})}{(F,L^{i_0})}}(r)+O(1)\leq T_F(r)+O(1)=T_f(r)+O(1).
\end{align*}
Therefore, by the Cartan-Nochka Second Main Theorem we have
\begin{align*}
T_f(r)+O(1)\geq \sum_{k=1,k\ne i_0}^qN_{(f,H_k)}^{[1]}(r)\geq \sum_{k=1,k\ne i_0}^q\frac{1}{p}N_{(f,H_k)}^{[p]}(r)\\
\geq \frac{q-2n+p-2}{p}T_f(r)-o(T_f(r)).
\end{align*}
This implies that $q\leq 2n+2.$ This is a contradiction. Hence, we get (\ref{a0}) in this case.

{\bf Case 2:} For any subset $J\subset \{1,\dots,q\}$ with $\#J= n+1,$ there exists a pair $i, j\in J$ such that $$\frac{(f,H_i)}{(g,L_i)}\not\equiv\frac{(f,H_j)}{(g,L_j)}.$$

\noindent We introduce an equivalence relation on $ L:=\{1,\cdots, q\}$ as follows: $i\sim j$ if and only if
\begin{equation*}
\text{det} 
\begin{pmatrix}
(f,H_{i}) &  & (f,H_{j})\cr&  \cr (g,L_{i}) &  & 
(g,L_{j})
\end{pmatrix}
\equiv 0.
\end{equation*}
Set $\{L_1,\cdots, L_s\}=L/\sim $. It is clear that  
 $\sharp L_k\leq n$ for all $k\in\{1,\cdots, s\}.$ Without loss of generality, we may assume that $L_k:=\{i_{k-1}+1,\cdots, i_k\}$
($k\in\{1,\cdots, s\})$ where $0=i_0<\cdots <i_s=q.$ 

\noindent We define the map $\sigma: \{1,\cdots, q\}\to \{1,\cdots, q\}$ by
\begin{equation*}
\sigma (i)=
\begin{cases}
i+n& \text{ if $i+n\leq q$},\\
i+n-q& \text{ if  $i+n>q$}.
\end{cases}
\end{equation*}
It is easy to see that  $\sigma$ is bijective and $\mid \sigma (i)-i\mid\geq n $ (note  that $q> 2n+2).$ This implies that $i$ and $\sigma (i)$ belong to distinct sets of $\{L_1,\cdots, L_s\}.$ This implies that for all $i\in\{1,\dots,q\},$
\begin{equation*}
P_i:=\text{det} 
\begin{pmatrix}
(f,H_{i}) &  & (f,H_{\sigma(i)})\cr&  \cr (g,L_{i}) &  & 
(g,L_{\sigma(i)})
\end{pmatrix}
\not\equiv 0.
\end{equation*}
By the assumption and by the definition of function $P_i,$  we have
\begin{align}\label{h2}
\nu_{P_i}\geq \min\{\nu_{(f,H_i)},\nu_{(g,L_i)}\}+\min\{\nu_{(f,H_{\sigma(i)})},\nu_{(g,L_{\sigma(i)})}\}+\sum_{j=1\atop j\ne i,\sigma(i)}^q\nu^{[1]}_{(f,H_j)}
\end{align}
outside an analytic set of codimension $\geq 2.$

\noindent On the other hand, since $f^{-1}(H_k)=g^{-1}(L_k)$ we have 
\begin{align*}
\min\{\nu_{(f,H_k)},\nu_{(g,L_k)}\}&\geq \min\{\nu_{(f,H_k)},p\} +\min\{\nu_{(g,L_k)},p\}-p\min\{\nu_{(f,H_k)},1\}\\
&=\nu^{[p]}_{(f,H_k)}+\nu^{[p]}_{(g,L_k)}-p\nu^{[1]}_{(f,H_k)}
\end{align*}
for $k\in\{i,\sigma(i)\}.$ 

\noindent Therefore, by (\ref{h2}) we have
\begin{align*}
\nu_{P_i}\geq \nu_{(f,H_i)}^{[p]}& +\nu_{(g,L_i)}^{[p]}+\nu_{(f,H_{\sigma(i)})}^{[p]}+\nu_{(g,L_{\sigma(i)})}^{[p]}\notag\\
&-p\nu_{(f,H_i)}^{[1]}-p\nu_{(f,H_{\sigma(i)})}^{[1]}+\sum_{j=1\atop j\ne i,\sigma(i)}^q\nu^{[1]}_{(f,H_j)}
\end{align*}
outside an analytic set of codimension $\geq 2.$

\noindent Then for all  $i\in\{1,\dots,q\}$ we have
\begin{align}\label{h3}
N_{P_i}(r)\geq N_{(f,H_i)}^{[p]}(r)& +N_{(g,L_i)}^{[p]}(r)+N_{(f,H_{\sigma(i)})}^{[p]}(r)+N_{(g,L_{\sigma(i)})}^{[p]}(r)\notag\\
&-pN_{(f,H_i)}^{[1]}(r)-p N_{(f,H_{\sigma(i)})}^{[1]}(r)+\sum_{j=1\atop j\ne i,\sigma(i)}^q N^{[1]}_{(f,H_j)}(r).
\end{align}
On the other hand, by Jensen's formula 
\begin{align*}
N_{P_i}(r)&=\int_{S(r)}\log |P_i| \sigma+O(1)\\
&\leq \int_{S(r)}\log(|(f,H_i)|^2+|(f,H_{\sigma(i)})|^2)^{\frac{1}{2}}\sigma\\
&\quad\quad\quad\quad\quad\quad\quad\quad\quad+\int_{S(r)}\log(|(g,L_i)|^2+|(g,L_{\sigma(i)})|^2)^{\frac{1}{2}}\sigma+O(1)\\
&\leq T_f(r)+T_g(r)+O(1).
\end{align*}
Therefore, by (\ref{h3}) for  all $i\in\{1,\dots,q\}$ we have
\begin{align}\label{h3'}
 N_{(f,H_i)}^{[p]}(r)& +N_{(g,L_i)}^{[p]}(r)+N_{(f,H_{\sigma(i)})}^{[p]}(r)+N_{(g,L_{\sigma(i)})}^{[p]}(r)\notag\\
&-pN_{(f,H_i)}^{[1]}(r)-p N_{(f,H_{\sigma(i)})}^{[1]}(r)+\sum_{j=1\atop j\ne i,\sigma(i)}^q N^{[1]}_{(f,H_j)}(r)\notag\\
&\leq T_f(r)+T_g(r)+O(1).
\end{align}
By summing-up of both sides of the above inequality for all $ i\in\{1,\dots,q\},$ we have
\begin{align}\label{h4}
2\sum_{j=1}^q\big(N_{(f,H_j)}^{[p]}(r)+N_{(g,L_j)}^{[p]}(r)\big)&+(q-2p-2)\sum_{j=1}^qN^{[1]}_{(f,H_j)}(r)\notag\\
&\leq q\big( T_f(r)+T_g(r)\big)+O(1).
\end{align}
Therefore, since $f^{-1}(H_j)=g^{-1}(L_j)$   we have
\begin{align}\label{h5}
2\sum_{j=1}^q\big(N_{(f,H_j)}^{[p]}(r)+N_{(g,L_j)}^{[p]}(r)\big)&+\frac{q-2p-2}{2}\sum_{j=1}^q\big(N^{[1]}_{(f,H_j)}(r)+N^{[1]}_{(g,L_j)}(r)\big)\notag\\
&\leq q\big( T_f(r)+T_g(r)\big)+O(1).
\end{align}
Then
\begin{align}\label{h6}
\big (2+\frac{q-2p-2}{2p}\big)\sum_{j=1}^q\big(N_{(f,H_j)}^{[p]}(r)+N_{(g,L_j)}^{[p]}(r)\big)
\leq q\big( T_f(r)+T_g(r)\big)+O(1).
\end{align}
 By (\ref{h6})  and  by  the Cartan-Nochka Second Main Theorem  we have
\begin{align*}
\frac{(q+2p-2)(q-2n+p-1)}{2p}\big(T_f(r)+T_g(r)\big)
&\leq q\big( T_f(r)+T_g(r)\big)+o\big( T_f(r)+T_g(r)\big).
\end{align*}
It follows that $(q+2p-2)(q-2n+p-1)\leq 2pq.$ Then $q^2-(2n+3-p)q-2(p-1)(2n+1-p)\leq 0.$ This is a contradiction to condition $(*)$ of Theorem \ref{TQT}. 
Thus we have completed the proof of (\ref{a0}).\\

 Assume that $H_j: a_{j0}\omega_0+\cdots+a_{jn}\omega_n=0,$  $L_j: b_{j0}\omega_0+\cdots+b_{jn}\omega_n=0\quad  (j=1,\dots,q).$ 

Set
$$A:=\begin{pmatrix}
a_{10}& \dots & a_{1n}\\
a_{20}& \dots & a_{2n}\\
\vdots & \ddots & \vdots\\
a_{(n+1)0}& \dots & a_{(n+1)n} \end{pmatrix},   B:=  \begin{pmatrix}
b_{10}& \dots & b_{1n}\\
b_{20}& \dots & b_{2n}\\
\vdots & \ddots & \vdots\\
b_{(n+1)0}& \dots & b_{(n+1)n} \end{pmatrix},\;\text{and}\;\mathcal L=B^{-1}\cdot A.$$

\noindent By (\ref{a0}), we have $A(f)\equiv B(g)$,  so we get  $\mathcal L(f)\equiv g.$

 Set $H^*_j=(a_{j0},\dots,a_{jn})\in \C^{n+1},$ $L^*_j=(b_{j0},\dots,b_{jn})\in \C^{n+1}.$ We write $H^*_j=\alpha_{j1}H^*_1+\cdots+\alpha_{j(n+1)}H^*_{n+1}$ and $L^*_j=\beta_{j1}L^*_1+\cdots+\beta_{j(n+1)}L^*_{n+1}.$

\noindent By (\ref{a0}) we have
\begin{align*}
\frac{\alpha_{j1}(f,H_1)+\dots+\alpha_{j(n+1)}(f,H_{n+1})}{\beta_{j1}(g,L_1)+\dots+\beta_{j(n+1)}(g,L_{n+1})}\equiv\frac{(f,H_1)}{(g,L_1)}\equiv\cdots\equiv\frac{(f,H_{n+1})}{(g,L_{n+1})}
\end{align*}
 for all $j\in\{1,\dots,q\}.$

\noindent  This implies that
\begin{align}\label{tp}
(\alpha_{j1}-\beta_{j1})(f,H_1)+\dots+(\alpha_{j(n+1)}-\beta_{j(n+1)})(f,H_{n+1})\equiv 0
\end{align}
for all $j\in\{1,\dots,q\}.$

\noindent On the other hand $f :\C^m\longrightarrow \langle\text{Im}f\rangle$ is linearly nondegenerate and $\{H_j\}_{j=1}^{n+1}$ are in general position  in $\C P^n.$  Thus, by (\ref{tp})  we have
\begin{align}\label{tp1}
 (\alpha_{j1}-\beta_{j1})(\omega, H_1)+\dots+(\alpha_{j(n+1)}-\beta_{j(n+1)})(\omega, H_{n+1})=0
\end{align}
 for all $\omega\in  \langle\text{Im}f\rangle$ for all $j\in\{1,\dots,q\}.$ 

Let hyperplanes $\alpha_j: \alpha_{j1}\omega_0+\cdots+\alpha_{j(n+1)}\omega_n=0$ and $\beta_j: \beta_{j1}\omega_0+\cdots+\beta_{j(n+1)}\omega_n=0\quad (j=1,\dots,q).$

\noindent By (\ref{tp1}) we have 
\begin{align}\label{tp2}
(A(\omega),\alpha_j)=( A(\omega),\beta_j)
\end{align}
for all $\omega\in  \langle\text{Im}f\rangle$ and $ j\in\{1,\dots,q\}.$ 

For any $j\in\{1,\dots,q\}$ and for any $\omega\in\langle\text{Im}f\rangle $ we have
\begin{align*}
(\omega,H_j)&=\alpha_{j1}(\omega,H_1)+\dots+\alpha_{j(n+1)}(\omega,H_{n+1})\\
&=( A(\omega),\alpha_j)\\
&\stackrel{(\ref{tp2})}{=}( A(\omega), \beta_j)\\
&=( B\cdot\mathcal L(\omega),\beta_j)\\
&=\beta_{j1}(\mathcal L(\omega),L_1)+\dots+\beta_{j(n+1)}(\mathcal L(\omega),L_{n+1})\\
&=(\mathcal L(\omega),L_j).
\end{align*}
This implies that  $\mathcal L(\langle\text{Im}f\rangle\cap H_j)=L_j\cap\mathcal L(\langle\text{Im}f\rangle)$  for all $j\in\{1,\dots,q\},$
which completes the proof of Theorem \ref{TQT}.
\hfill$\Box$

\vspace{1cm}

 \noindent  Gerd Dethloff$^{1-2} $ \\
 $^1$ Universit\'e Europ\'eenne de Bretagne, France\\
 $^2$
Universit\'{e} de Brest \\
   Laboratoire de math\'{e}matiques \\
UMR CNRS 6205\\
6, avenue Le Gorgeu, BP 452 \\
   29275 Brest Cedex, France \\
e-mail: gerd.dethloff@univ-brest.fr\\

\noindent Si Duc Quang and Tran Van Tan  \\
Department of Mathematics\\
  Hanoi National University of Education\\
 136-Xuan Thuy street, Cau Giay, Hanoi, Vietnam\\
e-mails: ducquang.s@gmail.com and tranvantanhn@yahoo.com  


\begin{thebibliography}{99}
\bibitem{CY}Z. Chen and Q. Yan, \textit{ Uniqueness theorem of meromorphic maps into $\Bbb P^N(\Bbb C)$ sharing $2N+3$ hyperplanes regardless of multiplicities,}  Internat. J. Math.  \textbf{20}  (2009),   717--726. 

\bibitem{DT} G. Dethloff and T. V. Tan,  \textit{Uniqueness theorems for
meromorphic mappings with few hyperplanes,}  Bulletin
des Sciences Math\'{e}matiques, \textbf{133} (2009), 501-514.

\bibitem{F1}  {\ } H. Fujimoto, \textit{The uniqueness problem of meromorphic
maps into the complex projective space}, Nagoya Math. J. \textbf{58} (1975),
1-23.

\bibitem{F2}  {\ } H. Fujimoto, \textit{Uniqueness problem with truncated
multiplicities in value distribution theory}, Nagoya Math. J. \textbf{152}
(1998), 131-152.

\bibitem{Ji}  {\ }  S. Ji, \textit{Uniqueness problem without multiplicities in
value distribution theory}, Pacific J. Math. \textbf{135} (1988), 323-348.

\bibitem{Ne}  {\ }  R. Nevanlinna, \textit{Einige Eindeutigkeitss\"{a}tze in der
Theorie der meromorphen Funktionen}, Acta. Math. \textbf{48} (1926), 367-391.

\bibitem{No}  {\ } J. Noguchi, \textit{A note on entire pseudo-holomorphic curves and the proof of Cartan-Nochka's theorem,} Kodai Math. J. \textbf{28} (2005), 336-346.


\bibitem{Q}{\ }  S. D. Quang, \textit{Unicity problem
of meromorphic mappings sharing few hyperplanes},  preprint.


\bibitem{Sm}  {\ }  L. Smiley, \textit{Geometric conditions for unicity of
holomorphic curves}, Contemp. Math. \textbf{25} (1983), 149-154.

\bibitem{St}  {\ } W. Stoll, \textit{On the propagation of dependences}, Pacific J. Math. \textbf{139} (1989), 311-337.

\bibitem{TQ}{\ }  D. D. Thai and S. D. Quang, \textit{Uniqueness problem with 
truncated multiplicities
of meromorphic mappings in several complex variables}, Inter. J. Math., \textbf{17} (2006), 1223-1257.



\end{thebibliography}
\end{document}